\documentclass[10pt]{article}
\usepackage{macros_colas2}

\newtheorem{theorem}{Theorem}[section] 

\newtheorem{definition}[theorem]{Definition}

\newtheorem{lemma}[theorem]{Lemma}
\newtheorem{fact}[theorem]{Fact}

\newcommand{\TGro}{T^{\mathrm{(Gro)}}}
\newcommand{\TZar}{T^{\mathrm{(Zar)}}}

\title{A sufficient condition for  the  fiber of the tangent bundle of a scheme and its Zariski tangent space to be isomorphic}

\author{%
Colas Bardavid \\
IMSc --- 
CIT Campus, Taramani \\
Chennai 600\,113 India
}

\begin{document}

\maketitle

 \begin{center}  
\rule{10 cm}{.5pt} \end{center}
\begin{small}
\begin{center}
\begin{minipage}{10cm}
\textbf{Abstract --} 
In this note, we give a simple sufficient condition for the Zariski relative tangent space and the Grothendieck relative tangent space to be isomorphic.
\end{minipage}
\end{center}
\end{small}
 \begin{center}  
\rule{10 cm}{.5pt} \end{center}

\begin{center}
 \begin{minipage}{10cm}
\end{minipage}
\end{center}

\newpage

Two notions of tangent space have been proposed in  scheme theory: the Zariski tangent space $T_{x}X$ and the Grothendieck relative tangent space $\TGro_{X/S}(x)$. The relation between these two tangent spaces is known only in a very special case: in \cite{EGA44}, Grothendieck shows that for a sheme $X$, when $S=\Sp k$ (where $k$ is a field) and when $x\in X$ is $k$-rational point of $X$, then these two objects coincide.  The aim of this note is to compare these two tangent spaces in more general situations. After having  introduced a relative tangent space in the Zariski fashion,  we show that there always exists a morphism from the Grothendieck relative tangent space to the Zariski one. The main result is that this morphism is an isomorphism whenever the extension $\kappa(x)/\kappa(s)$ is algebraic and separable. We also give a counter-example showing that in general these two tangent spaces do not coincide.

\medskip

The note is organized as follows: in Section 1, we introduce three different tangent spaces $T_{x}X$, $\TGro_{X/S}(x)$ and $\TZar_{X/S}(x)$. In Section 2, we prove that in general $\TGro_{X/S}(x)$ and $\TZar_{X/S}(x)$ are not isomorphic. In section 3, we construct a morphism $\TGro_{X/S}(x) \longto \TZar_{X/S}(x)$. In section 4, we give  a condition for this morphism to be an isomorphism. In Section 5, we prove the main result of this note: $\TGro_{X/S}(x) $ and $ \TZar_{X/S}(x)$ are isomorphic whenever $\kappa(x)/\kappa(s)$ is an algebraic and separable extension. 

\medskip

We wish to thank David Madore for a helpful discussion on this matter.

\section{Introducing the two tangent spaces}

\subsection{Global notations. The Zariski tangent space. }In this paragraph, we recall some very classical facts and set various notations. Let $X\longto S$ be  schemes and $x\in X$ an element above $s\in S$. The structure sheaf of $X$ is denoted by $\mathcal{O}_{X}$, its stalk over $x$ by $\mathcal{O}_{X,x}$. Its maximal ideal is denoted by $\mathfrak{M}_{x}$ and  $\kappa(x):=\mathcal{O}_{X,x}/\mathfrak{M}_{x}$. For any $f\in \mathcal{O}_{X,x}$, the image of $f$ in $\kappa(x)$ under the canonical projection is denoted by $f(x)$. The residue field $\kappa(x)$ is viewed as a $\mathcal{O}_{X,x}$-module \emph{via}
$f\cdot \lambda:=f(x)\lambda$. The ideal $\mathfrak{M}_{x}$ is an $\mathcal{O}_{X,x}$-module and admits as a sub-$\mathcal{O}_{X,x}$ -module the ideal ${\mathfrak{M}_{x}}^2$. The quotient $\mathcal{O}_{X,x}$-module $\mathfrak{M}_{x}/{\mathfrak{M}_{x}}^2$ is in fact a $\kappa(x)$-vector space. We will also denote by $[f]$ the image of $f\in \mathfrak{M}_{x}$ in $\mathfrak{M}_{x}/{\mathfrak{M}_{x}}^2$.

\begin{definition}
The \emph{Zariski tangent space of $X$ at $x$} is the $\kappa(x)$-vector space
$$
T_{x}X:=\Hom_{\kappa(x)}(\mathfrak{M}_{x}/{\mathfrak{M}_{x}}^2, \kappa(x)).
$$
\end{definition}

\subsection{The Grothendieck relative tangent space. }It is defined in \S \textbf{(16.5.13)} of \cite{EGA44} as follows.

\begin{definition}
 The \emph{Grothendieck relative tangent space of $X/S$ at $x$} is the $\kappa(x)$-vector space
 $$
 \TGro_{X/S}(x):=\Hom_{\kappa(x)}  ( \Omega_{X/S}^1\otimes_{\mathcal{O}_X}\kappa(x), \kappa(x)  ).
 $$
\end{definition}

Let us recall that $\Omega^1_{X/S}$ is the $\mathcal{O}_X$-module of $1$-differentials of $X/S$. We don't go further in the description of this object since Fact \ref{eq.def.Gro.tg} will give a more handy definition of $\TGro_{X/S}(x)$.

The main advantage of this construction over the Zarsiki's one is that $ \TGro_{X/S}(x)$ appears as the fiber of a tangent bundle. The \emph{tangent bundle of $X$ relatively to $S$} is  a vector bundle $T_{X/S}$ above $X$. By definition, the $\kappa(x)$-rational points of its fiber over $x$ is precisely $ \TGro_{X/S}(x)$. Let us give now the more handy  description of $ \TGro_{X/S}(x)$, which has also been noticed by Kunz in \cite{Onthetangentbundle}:

\begin{fact}\label{eq.def.Gro.tg}$ \TGro_{X/S}(x)\simeq \mathrm{Der}_{\mathcal{O}_{S,s}} \big ( \mathcal{O}_{X,x},\, \kappa(x) \big ).$\label{new.expression.Grot}
\end{fact}

\pproof{
It follows from two observations. First, as in \cite{EGA44}, one can writes
$$
\TGro_{X/S}(x)=\Hom_{\kappa(x)}\left ( \quotient{\Omega_{\mathcal{O}_{X,x}/\mathcal{O}_{S,s}}}{\mathfrak{M}_{x}\cdot \Omega_{\mathcal{O}_{X,x}/\mathcal{O}_{S,s}}}\ , \kappa(x)\right ).
$$
Second, if $k$ is a ring and $(A,\mathfrak{M}) $ a local $k$-algebra with  residual field $K$, then 
$$\Hom_{K-ev}\left ( \Omega_{A/k} /\mathfrak{M}\cdot \Omega_{A/k}, K \right )\simeq \mathrm{Der}_{k}\left ( A, K \right ),$$
as one can easily check. Then, apply this formula with $k=\mathcal{O}_{S,s}$, $A=\mathcal{O}_{X,x}$ to conclude. }

\subsection{The Zariski relative tangent space. } Let us give a definition of the relative tangent space, in the Zariski fashion. When $X$ and $S$ are schemes, with $f : X\longto S$ sending $x$ to $s$, one would like to define the differential of $f$ in $x$, mapping $T_{x}X$ to $T_{s}S$ (or better, to $T_{s}S\otimes_{\kappa(s)}\kappa(x)$). Imagine there is such a $\kappa(x)$-linear map
$$
T_{x}f : T_{x}X \longto T_{s}S\otimes_{\kappa(s)}\kappa(x).
$$
Then, the relative tangent space would be the kernel of this map. Intuitively, it corresponds to the tangent space of the fiber $X_{s}$ at $x$. Actually, such a map $T_{x}f$ \emph{does not} exist, but we can still define a very similar morphism and, subsequently, the relative tangent space.
\medskip

The morphism $f$ induces a morphism $f_{x}^\# : \mathcal{O}_{S,s}\longto \mathcal{O}_{X,x}$
that maps $\mathfrak{M}_{s}$ into $\mathfrak{M}_{x}$. Hence, we get a map
$
 i_{x} : \kappa(s)\longto \kappa(x)
$ which is an injection between fields, as well as
$$
j_{x} : \mathfrak{M}_{s}/{\mathfrak{M}_{s}}^2\longto \mathfrak{M}_{x}/{\mathfrak{M}_{x}}^2
$$
which is $\kappa(s)$-linear. 
So, starting with $\vec{v}\in T_{x}X$, one obtains a $\kappa(s)$-linear map that we will denote
$$
\widetilde{T}_{x}f \bullet \vec{v} :=\vec{v}\circ j_{x}\in \Hom_{\kappa(s)}\left ( \mathfrak{M}_{s}/{\mathfrak{M}_{s}}^2, \kappa(x) \right ).
$$
The application
$$
\widetilde{T}_{x}f : \fonctionb{T_{x}X}{\Hom_{\kappa(s)}(\mathfrak{M}_{s}/{\mathfrak{M}_{s}}^2, \kappa(x))}{\vec{v}}{T_{x}f \bullet \vec{v}}
$$
is $\kappa(x)$-linear. The problem to define properly the differential of $f$ is that, in general, $\Hom_{\kappa(s)}(\mathfrak{M}_{s}/{\mathfrak{M}_{s}}^2, \kappa(x))$ and $T_{s}S \otimes_{\kappa(s)}\kappa(x)$ are not isomorphic. But here, the map $\widetilde{T}_{x}f$ will play the role of the differential.

\smallskip

\begin{definition}
The \emph{Zariski relative tangent space of $X/S$ at $x$} is the $\kappa(x)$-vector space
$$
\TZar_{X/S}(x):= \ker \, \widetilde{T}_{x}f.
$$
\end{definition}

\subsection{An alternative description of the Zariski relative tangent space.} As intuition suggests, the relative Zariski tangent space can be described as in the following lemma. We will use later this description.

\begin{lemma}\label{relativetoabsolute}
$\TZar_{X/S}(x)$ is isomorphic to the tangent space $T_{X_s,x}$ at $x$ of the fiber $X_s$ above $s$.
\end{lemma}

\medskip

\pproof{
First, let us describe $T_{X_s,x}$. The local ring $\mathcal{O}_{X_s, x}$ is isomorphic to $\mathcal{O}_{X,x} / \mathfrak{M}_s \mathcal{O}_{X,x}$. Its maximal ideal is $\mathfrak{M}_x / \mathfrak{M}_s \mathcal{O}_{X,x}$ and the tangent space is
$$
T_{X_s,x} \simeq \Hom_{\kappa(x)} \Big (\frac{\mathfrak{M}_x / \mathfrak{M}_s \mathcal{O}_{X,x}}{{\mathfrak{M}_x}^2 / \mathfrak{M}_s \mathcal{O}_{X,x}} , \kappa(x) \Big ).
$$
As just seen, the relative Zariski tangent space can be described as
$$
\TZar_{X/S}(x) = \Hom_{\kappa(x)} \Big ( \frac{\mathfrak{M}_x /{\mathfrak{M}_x}^2}{\mathfrak{M}_s/ {\mathfrak{M}_s}^2}, \kappa(x) \Big ).
$$
So, let us prove that 
$$
\frac{\mathfrak{M}_x / \mathfrak{M}_s \mathcal{O}_{X,x}}{{\mathfrak{M}_x}^2 / \mathfrak{M}_s \mathcal{O}_{X,x}} \simeq \frac{\mathfrak{M}_x /{\mathfrak{M}_x}^2}{\mathfrak{M}_s/ {\mathfrak{M}_s}^2}.
$$
Clearly, there is a map
$$
\mathfrak{M}_x/{\mathfrak{M}_x}^2 \longto \frac{\mathfrak{M}_x / \mathfrak{M}_s \mathcal{O}_{X,x}}{{\mathfrak{M}_x}^2 / \mathfrak{M}_s \mathcal{O}_{X,x}},
$$
whose kernel can be described as the image of $\mathfrak{M}_s \longto \mathfrak{M}_x/{\mathfrak{M}_x}^2$. To conclude, remark that this last map factors through $\mathfrak{M}_s/ {\mathfrak{M}_s}^2$. }

\section{Grothendieck and Zariski are not isomorphic in general}

At this point, it is very easy to give a counter-example. Indeed, let us consider $k$ a field and $X:=\Sp k(t)$, $S:=\Sp k$ and $x$ the unique element of $X$. Then, $T_x X=0$ and so $\TZar_{X/S}(x)=0$. But, 
$$
\TGro_{X/S}(x)=\mathrm{Der}_{k}(k(t), k(t))
$$
which is isomorphic to $k(t)$.

\section{From Grothendieck to Zariski}

 Let $D\in \TGro_{X/S}(x)$, in other words, in virtue of Fact \ref{new.expression.Grot}, let $D : \mathcal{O}_{X,x}\longto \kappa(x)$ be a $\mathcal{O}_{S,s}$-derivation.
We can associate to $D$ an element of $\TZar_{X/S}(x)$. Indeed, the restriction of $D$ to $\mathfrak{M}_{x}$ factors through $\mathfrak{M}_{x}\longto \mathfrak{M}_{x}/{\mathfrak{M}_{x}}^2$, since for any $f,g\in \mathfrak{M}_{x}$ one has
\begin{align*}
D(fg)=f\cdot D(g) + g\cdot D(f)  = f(x) \cdot D (g) + g(x) \cdot D(f)=0.
\end{align*} If  $$\phi_{D} : \fonctionb{\mathfrak{M}_{x}/{\mathfrak{M}_{x}}^2}{\kappa(x)}{[\varphi]}{D(\varphi)}$$ denotes the factored map, then let us check that $\phi_D \in \TZar_{X/S}(x)$. We have to show that the morphism
$$
\xymatrix{
{\mathfrak{M}_{s}/{\mathfrak{M}_{s}}^2} \ar[r]^-{j_x} & {\mathfrak{M}_{x}/{\mathfrak{M}_{x}}^2} \ar[r]^-{\phi_D} & \kappa(x)
}
$$
is zero. This is straighforward since $D$ is zero on $\mathcal{O}_{S,s}$.
Hence, we have defined a $\kappa(x)$-linear map
$$
\boldsymbol{\Phi}_{X/S}^x := \fonctionb{\TGro_{X/S}(x)}{\TZar_{X/S}(x)}{D}{\phi_{D}}.
$$

\section{A condition for Grothendieck and Zariski to be  isomorphic}

To begin with, let us construct analogs of $\mathfrak{M}_{x}$ and $\mathfrak{M}_{x}/{\mathfrak{M}_{x}}^2$ with a structure of  $\kappa(x)$-algebra. First, let us denote
$$
\widetilde{\mathcal{O}}_{X/S,\,x}:=\kappa(x)\otimes _{\mathcal{O}_{S,s} }\mathcal{O}_{X,x}
$$
and point out some facts:
\begin{itemize}
\item Any derivation $D\in \mathrm{Der}_{\mathcal{O}_{S,s}} \big ( \mathcal{O}_{X,x},\, \kappa(x) \big )$ gives rise to a derivation 
$$\widetilde{D} \in \mathrm{Der}_{\kappa(x)} \big ( \widetilde{\mathcal{O}}_{X/S,x},\, \kappa(x) \big ).$$
It is simply defined by $\widetilde{D}(\lambda\otimes \varphi) = \lambda\otimes D(\varphi)$.
\item On the ring $\widetilde{\mathcal{O}}_{X/S,\,x}$, we still have an evaluation map. It is
$$
\fonctionb{\widetilde{\mathcal{O}}_{X/S,\,x}}{\kappa(x)}{\lambda\otimes \varphi}{\lambda \cdot \varphi(x)}.
$$
 We then define
$
\widetilde{\mathfrak{M}}_{x/s}:= \ker (\widetilde{\mathcal{O}}_{X/S,\,x}\longto \kappa(x)).
$
Remark that as in classical case, we have $(\widetilde{\mathcal{O}}_{X/S,\,x}) / (\widetilde{\mathfrak{M}}_{x/s}) \simeq \kappa(x)$, so that the $\widetilde{\mathcal{O}}_{X/S,\,x}$-module $\widetilde{\mathfrak{M}}_{x/s}/{(\widetilde{\mathfrak{M}}_{x/s})}^2$ is actually a $\kappa(x)$-vector space.

\item If  $\ell_{x} : {\mathcal{O}_{X,x}}\longto {\widetilde{\mathcal{O}}_{X/S,\,x}}$  denotes the ring morphism that sends $\varphi$ to $1\otimes \varphi$, then $\ell_{x}$ maps $\mathfrak{M}_{x}$ into $\widetilde{\mathfrak{M}}_{x/s}$ an so, one can consider the following
$$
\vartheta_{x/s} :\mathfrak{M}_{x}/{\mathfrak{M}_{x}}^2\longto \widetilde{\mathfrak{M}}_{x/s}/{(\widetilde{\mathfrak{M}}_{x/s})}^2,
$$
which is a morphism of $\kappa(x)$-vector spaces. It makes the following diagram commute : 
\begin{equation}\label{Diag.Comm}
\xymatrix@C=1.3cm{
\ar[d]\mathfrak{M}_{x}\ar[r]^-{\ell_{x}}& \widetilde{\mathfrak{M}}_{x/s}\ar[d] \\
\mathfrak{M}_{x}/ {\mathfrak{M}_{x}}^2 \ar[r]^-{\vartheta_{x/s}} & \widetilde{\mathfrak{M}}_{x/s}/ ({\widetilde{\mathfrak{M}}_{x/s}})^2
}
\end{equation}
\item One has $\widetilde{D}(\ell_{x}(\varphi))=D(\varphi)$ and, as it happens for $D$, one can factor $\widetilde{D}$ through $\widetilde{\mathfrak{M}}_{x/s}/({\widetilde{\mathfrak{M}}_{x/s}})^2$.
\item  We still denote by $[\varphi]$ the image of $\varphi\in \widetilde{\mathfrak{M}}_{x/s}$ in $\widetilde{\mathfrak{M}}_{x/s}/{(\widetilde{\mathfrak{M}}_{x/s})}^2$.
\end{itemize} 
\bigskip

\noindent Now, let us state and prove the main lemma of this note.

\begin{lemma}
If 
 $\vartheta_{x/s}$ is an  isomorphism, then $\boldsymbol{\Phi}_{X/S}^x$ is also one.
\end{lemma}
\smallskip

\pproof{Let us assume that $\vartheta_{x/s}$ is an isomorphism. First, we construct a $\kappa(x)$-linear map $\boldsymbol{\Upsilon}^{x}_{X/S}:\TZar_{X/S}(x)\longto \TGro_{X/S}(x)$. Let $\vec{v}\in \TZar_{X/S}(x)$. We associate to $\vec{v}$ the following map
$$
D_{\vec{v}} : \fonctionb{\mathcal{O}_{X, x}}{\kappa(x)}
{\varphi}{\vec{v}\bullet {(\vartheta_{x/s})}^{-1}[1\otimes \varphi-\varphi(x)\otimes 1]}.
$$
It is an $\mathcal{O}_{S,s}$-derivation. Indeed,\medskip
\begin{itemize}
\item First, if $\varphi\in \mathcal{O}_{S,s}$ then $1\otimes \varphi=\varphi(x)\otimes 1$ and so
\begin{align*}
D_{\vec{v}} (f_{x}^\#(\varphi))&= \vec{v}\bullet {(\vartheta_{x/s})}^{-1}[1\otimes \varphi-\varphi(x)\otimes 1]\\
&= \vec{v}\bullet {(\vartheta_{x/s})}^{-1}[\varphi(x)\otimes 1-\varphi(x)\otimes 1]=0.
\end{align*}
\item  Second, let us verify the Leibniz rule. Let $\varphi, \psi\in \mathcal{O}_{X,x}$. In $\widetilde{\mathfrak{M}}_{x/s}/{(\widetilde{\mathfrak{M}}_{x/s})}^2$, one has
$$
[1\otimes \varphi-\varphi(x)\otimes 1]\cdot [1\otimes\psi-\psi(x)\otimes 1]=0
$$
and so
\begin{equation}\label{une.equation}
[\psi(x)\otimes \varphi]+ [\varphi(x)\otimes \psi]
 =[1\otimes \varphi \psi] + [\varphi(x)\psi(x)\otimes 1] 
\end{equation}
Then, 
\begin{align*}
&\psi(x) \cdot [1\otimes \varphi-\varphi(x)\otimes 1] + 
\varphi(x) \cdot [1\otimes \psi-\psi(x)\otimes 1] \\
&= [\psi(x)\otimes \varphi] + [\varphi(x)\otimes \psi] - 2 [\varphi(x)\psi(x)\otimes 1] \\
&=[1\otimes \varphi \psi] - [\varphi(x)\psi(x)\otimes 1]\qquad \qquad \qquad \qquad \qquad  \text{by (\ref{une.equation}).}
\end{align*}
This implies that $D_{\vec{v}}(\varphi \psi)=\varphi(x) \cdot D_{\vec{v}}(\psi) + \psi(x) \cdot D_{\vec{v}}(\varphi)$.
\end{itemize}
\smallskip

Now, let us prove that $\boldsymbol{\Upsilon}^x_{X/S}\circ \boldsymbol{\Phi}_{X/S}^x=\Id$. Let $D : \mathcal{O}_{X,x}\longto \kappa(x)$ be an $\mathcal{O}_{S,s}$-derivation. Let $\varphi\in \mathcal{O}_{X,x}$. Let us compute 
\begin{align*}
\boldsymbol{\Upsilon}^x_{X/S}\circ \boldsymbol{\Phi}_{X/S}^x (D)(\varphi) &= \boldsymbol{\Phi}_{X/S}^x (D) \bullet  (\vartheta_{x/s})^{-1}  [1\otimes \varphi - \varphi(x)\otimes 1].
\end{align*}
Let $\psi\in \mathcal{O}_{X,x}$ such that
\begin{equation}\label{Defi.Psi}
[\psi]=({\vartheta_{x/s}})^{-1} [1\otimes \varphi - \varphi(x)\otimes 1].
\end{equation}
Then, by definition, one has
$ \boldsymbol{\Phi}_{X/S}^x (D) \bullet (\vartheta_{x/s})^{-1} [1\otimes \varphi - \varphi(x)\otimes 1] = D(\psi)$. Applying $\vartheta_{x/s}$ to (\ref{Defi.Psi}), one gets, with (\ref{Diag.Comm}), that 
$$
[1\otimes \varphi-\varphi(x)\otimes 1] = [1\otimes \psi].
$$
Applying $\widetilde{D}$, and since  $\widetilde{D}(\varphi(x)\otimes 1)=0$, one obtains that
$$
\widetilde{D} ([1\otimes \varphi-\varphi(x)\otimes 1]) = D(\varphi) = D(\psi).
$$
So, we have got the required identity, $\boldsymbol{\Upsilon}^x_{X/S}\circ \boldsymbol{\Phi}_{X/S}^x (D)=D$.

\medskip

Let us prove now that $\boldsymbol{\Phi}_{X/S}^x \circ \boldsymbol{\Upsilon}^x_{X/S} =\Id$. Let $\vec{v} : \mathfrak{M}_x / {\mathfrak{M}_x }^2 \longto \kappa(x)$ be a Zariski tangent vector and let $\varphi\in \mathfrak{M}_x$. Then, 
\begin{align*}
  \boldsymbol{\Phi}_{X/S}^x  \circ \boldsymbol{\Upsilon}^x_{X/S}(\vec{v})\bullet \varphi & =  \boldsymbol{\Upsilon}^x_{X/S} (\vec{v}) (\varphi) \\
  &= \vec{v}\bullet (\vartheta_{x/s})^{-1}  [1\otimes \varphi - \varphi(x)\otimes 1] \\
  &= \vec{v}\bullet (\vartheta_{x/s})^{-1}  [1\otimes \varphi] =  \vec{v}\bullet ((\vartheta_{x/s})^{-1} \circ \vartheta_{x/s} )( \varphi) \\ &= \vec{v}\bullet \varphi
  \end{align*}}

\section{The main theorem}

\begin{theorem}
  When the extension $i_x : \kappa(s)\longto \kappa(x)$ is algebraic and separable
$$
\boldsymbol{\Phi}_{X/S}^x : \xymatrix@1{
\TGro_{X/S}(x) \ar[r]^-{\widetilde{\phantom{aaa}}} & \TZar_{X/S}(x)
}
$$
is an isomorphism of $\kappa(x)$-vector spaces.
\end{theorem}

\medskip

\pproof{To begin with, remark that we can replace the relative situation $X\longto S$ by the ``absolute situation'' $X_s \longto \Sp \kappa(s)$. Indeed, let us consider the following cartesian square
$$
\xymatrix@R=5mm@C=4mm{
X_s \ar[d] \ar[r] & X \ar[d] \\
\Sp \kappa(s) \ar[r] & S
}.
$$
First, by  Lemma \ref{relativetoabsolute}, the relative Zariski tangent spaces are  isomorphic one to each other in this case. For the  Grothendieck tangent spaces, one can say  the following:
\begin{itemize}\setlength{\itemsep}{0mm}
\item[---] By \textbf{(16.5.13.2)} of \cite{EGA44}, when $\Omega_{X/S}^1$ is an $\mathcal{O}_X$-module of finite type, the Grothendieck tangent space is invariant under base extension, and so $\TGro_{X/S}(x)$ and $\TGro_{X_s/\kappa(s)}(x)$ are isomorphic.
\item[---] But, actually, the latter is true without any condition of finiteness, as we prove it in Lemma \ref{lemm.final}.
\end{itemize}

\noindent So, in what follows, we will assume that $\mathcal{O}_{S,s}=\kappa(s)$. In particular, $\mathcal{O}_{X,x}$ is a $\kappa(s)$-algebra !
\medskip

Now, let us apply the second  fundamental exact sequence of K\"ahler differentials (Theorem 58 of \cite{MatsumuraCommAlg}), respectively with 
\begin{itemize}\setlength{\itemsep}{-1mm}
 \item[1)] first, $k=\kappa(s)$, $A=\mathcal{O}_{X,x}$ and $\mathfrak{M}=\mathfrak{M}_x$ 
 \item[2)] second, $k=\kappa(x)$, $A={\widetilde{\mathcal{O}}_{X/S,\,x}} $ and $\mathfrak{M}= \widetilde{\mathfrak{M}}_{x/s}$
\end{itemize}
 to get the following two exact sequences:
\begin{align}
 &\mathfrak{M}_x/ {\mathfrak{M}_x}^2 \longto \Omega_{\mathcal{O}_{X,x} / \kappa(s)} \otimes_{\mathcal{O}_{X,x}} \kappa(x) \longto \Omega_{\kappa(x)/ \kappa(s)} \longto 0 \label{ExSeq1}\\
 \text{and} \quad & \widetilde{\mathfrak{M}}_{x/s}/ (\widetilde{\mathfrak{M}}_{x/s})^2 \longto \Omega_{ {\widetilde{\mathcal{O}}_{X/S,\,x}} / \kappa(x)}\otimes_{ {\widetilde{\mathcal{O}}_{X/S,\,x}}} \kappa(x) \longto \Omega_{\kappa(x)/\kappa(x)}=0. \label{ExSeq2}
\end{align}
In the second one, the left-hand morphism is injective. Indeed, ${\widetilde{\mathcal{O}}_{X/S,\,x}}  \longto \kappa(x)$ has a section and so, the criterion given by Proposition 16.12 of \cite{Eisenbud} applies. So, (\ref{ExSeq2}) can be written
$$
 0 \longto \widetilde{\mathfrak{M}}_{x/s}/ (\widetilde{\mathfrak{M}}_{x/s})^2 \longto \Omega_{ {\widetilde{\mathcal{O}}_{X/S,\,x}} / \kappa(x)}\otimes_{ {\widetilde{\mathcal{O}}_{X/S,\,x}}} \kappa(x) \longto 0
$$
and hence gives an isomorphism. But we also know that 
$$
\Omega_{ {\widetilde{\mathcal{O}}_{X/S,\,x}} / \kappa(x)} \simeq  \Omega_{\mathcal{O}_{X,x}/\kappa(s)}\otimes_{\mathcal{O}_{X,x}} {\widetilde{\mathcal{O}}_{X/S,\,x}}
$$
so that, we have
$$
  \Omega_{ {\widetilde{\mathcal{O}}_{X/S,\,x}} / \kappa(x)}\otimes_{ {\widetilde{\mathcal{O}}_{X/S,\,x}}} \kappa(x) 
  \simeq
   \Omega_{\mathcal{O}_{X,x}/\kappa(s)}\otimes_{\mathcal{O}_{X,x}} \kappa(x).
$$
Hence, in the exact sequence (\ref{ExSeq1}), we can replace the second vector space by $ \widetilde{\mathfrak{M}}_{x/s}/ (\widetilde{\mathfrak{M}}_{x/s})^2$. We get 
$$
 \mathfrak{M}_x/ {\mathfrak{M}_x}^2 \longto \widetilde{\mathfrak{M}}_{x/s}/ (\widetilde{\mathfrak{M}}_{x/s})^2 \longto \Omega_{\kappa(x)/ \kappa(s)} \longto 0
$$
and one can check that the first arrow in this sequence is $\vartheta_{x/s}$. 

\smallskip

By Corollary 16.13 of \cite{Eisenbud}, a sufficient condition for $\vartheta_{x/s}$ to be injective is that the extension $\kappa(x)/\kappa(s)$ is separable. A sufficient condition for $\vartheta_{x/s}$ to be surjective is that $\Omega_{\kappa(x)/\kappa(s)}=0$. Hence, if $\kappa(x)/\kappa(s)$ is separable and algebraic, by Lemma 16.15 of \cite{Eisenbud} $\vartheta_{x/s}$ is an isomorphism and so is $ \boldsymbol{\Phi}_{X/S}^x$. }

\bigskip
\begin{lemma}\label{lemm.final}
For any schemes $X/S$ and any $x\in X$ above $s\in S$,
$$\TGro_{X/S}(x)\simeq \TGro_{X_s/\Sp \kappa(s)}(x).$$
\end{lemma}

\medskip
\pproof{
We use the same description of the local ring of $X_s$ at $x$ as in Lemma \ref{relativetoabsolute}. So, we want to prove that 
\[
\mathrm{Der}_{\mathcal{O}_{S,s}} (\mathcal{O}_{X,x}, \kappa(x))\qquad \text{and}\qquad
\mathrm{Der}_{\kappa(s)} (\mathcal{O}_{X,x}/\mathfrak{M}_s\mathcal{O}_{X,x}, \kappa(x))
\]
are isomorphic. The two inverse map are the most natural ones to describe, and the check that it works is left to the reader as an exercise.
}

\def\cprime{$'$}

\end{document}